\documentclass{article}
\usepackage{amsfonts}

\input{tcilatex}
\begin{document}

\begin{center}
\textbf{SPECIAL CURVES OF 4D GALILEAN SPACE}

\bigskip Mehmet BEKTA\c{S} Mahmut ERG\"{U}T Alper Osman \"{O}\u{G}RENM\.{I}%
\c{S}

F\i rat University, Faculty of Science, Department of Mathematics

23119 Elaz\i \u{g} / T\"{U}RK\.{I}YE

mbektas@firat.edu.tr mergut@firat.edu.tr ogrenmisalper@gmail.com

\bigskip
\end{center}

\textbf{Abstract.}{\small \ }Special curves and their characterizations are
one of the main area of mathematicians and physicians.

As a special curve we will mainly focus on Mannheim curve which has the
following relation:%
\[
k_{1}=\beta (k_{1}^{2}+k_{2}^{2}) 
\]%
where $k_{1}$ and $k_{2}$ are curvature and torsion, respectively.

In the present paper we define Mannheim curves for 4-dimensional Galilean
space and investigate some characterization of it.

\textbf{Keywords:}{\small \ Galilean Space, Mannheim Curve, Frenet Formula.}

\textbf{MSC 2000 classification:}{\small \ 53B30, 53A35.}

\bigskip

\textbf{2. PRELIMINARIES}\bigskip

The Galilean space is a 3D complex projective space $P_{3}$ in which the
absolute figure $\{w,f,I_{1},I_{2}\}$ consists of a real plane $w$ (the
absolute plane), a real line $f\subset w$ (the absolute line) and two
complex conjugate points $I_{1},I_{2}\in f$ (the absolute points).

The study of mechanics of plane-parallel motions reduces to the study of a
geometry of 3D space with coordinates $\{x,y,t\}$ are given by the motion
formula [8]. This geometry is called 3D Galilean geometry. In [8], is
explained that 4D Galilean geometry, which studies all properties invariant
under motions of objects in space, is even more complex.

In addition it is stated that this geometry can be described more precisely
as the study of those properties of 4D space with coordinates which are
invariant under the general Galilean transformations as follows:

\begin{eqnarray*}
x^{^{\prime }} &=&(\cos \beta \cos \alpha -\cos \gamma \sin \beta \sin
\alpha )x+(\sin \beta \cos \alpha -\cos \gamma \cos \beta \sin \alpha )y \\
&&+(\sin \gamma \sin \alpha )z+(v\cos \delta _{1})t+a \\
y^{^{\prime }} &=&-(\cos \beta \sin \alpha +\cos \gamma \sin \beta \cos
\alpha )x+(-\sin \beta \sin \alpha +\cos \gamma \cos \beta \cos \alpha )y \\
&&+(\sin \gamma \cos \alpha )z+(v\cos \delta _{2})t+b \\
z^{^{\prime }} &=&(\sin \gamma \sin \beta )x-(\sin \gamma \cos \beta
)y+(\cos \gamma )z+(v\cos \delta _{3})t+c \\
t^{^{\prime }} &=&t+d
\end{eqnarray*}%
with $\cos ^{2}\delta _{1}+\cos ^{2}\delta _{2}+\cos ^{2}\delta _{3}=1.$

Some fundamental properties of curves in 4D Galilean space, is given for the
purpose of the requirements in the next section

A curve in $G_{4}$ $(I\subset \mathbb{R}\longrightarrow G_{4})$ is given as
follows

\[
\alpha (t)=(x(t),y(t),z(t),w(t)), 
\]%
where $x(t),y(t),z(t),w(t)\in C^{4}$ (smooth functions) and $t\in I$. Let $%
\alpha $ be a curve in $G_{4},$ which is parameterized by arclength $t=s,$
and given in the following coordinate form

\[
\alpha (s)=(s,y(s),z(s),w(s)). 
\]

In affine coordinates the Galilean scalar product between two points $%
P_{i}=(x_{i1},x_{i2},x_{i3},x_{i4}),$ $i=1,2$ is defined by

\[
g(P_{1},P_{2})=\{%
\begin{array}{c}
\left\vert x_{21}-x_{11}\right\vert ,\text{ if }x_{11}\neq x_{21}, \\ 
\sqrt{(x_{22}-x_{12})^{2}+(x_{23}-x_{13})^{2}+(x_{24}-x_{14})^{2}}\text{ if }%
x_{11}=x_{21}.%
\end{array}%
\]

For the vectors $a=(a_{1},a_{2},a_{3},a_{4},),$ $%
b=(b_{1},b_{2},b_{3},b_{4},)\ $and $c=(c_{1},c_{2},c_{3},c_{4},),$ Galilean
cross product in $G_{4}$ is defined as follows:

\[
a\wedge b\wedge c=\left\vert 
\begin{array}{cccc}
0 & e_{2} & e_{3} & e_{4} \\ 
a_{1} & a_{2} & a_{3} & a_{4} \\ 
b_{1} & b_{2} & b_{3} & b_{4} \\ 
c_{1} & c_{2} & c_{3} & c_{4}%
\end{array}%
\right\vert 
\]%
where $e_{i}$ are the standard basis vectors.

In this paper, we denote the inner product of two vectors $a,b$ in the sense
of Galilean by the notation $<a,b>_{G}.$

Let $\alpha (s)=(s,y(s),z(s),w(s))$ be a curve parameterized by arclength $s$
in $G_{4}.$ For a $\alpha $ Frenet curve, the Frenet formulas can be given
as following form

\[
\left[ 
\begin{array}{c}
\mathbf{t}^{^{\prime }} \\ 
\mathbf{n}^{^{\prime }} \\ 
\mathbf{b}^{^{\prime }} \\ 
\mathbf{e}^{^{\prime }}%
\end{array}%
\right] =\left[ 
\begin{array}{cccc}
o & \kappa & 0 & 0 \\ 
0 & 0 & \tau & 0 \\ 
0 & -\tau & 0 & \sigma \\ 
0 & 0 & -\sigma & 0%
\end{array}%
\right] \left[ 
\begin{array}{c}
\mathbf{t} \\ 
\mathbf{n} \\ 
\mathbf{b} \\ 
\mathbf{e}%
\end{array}%
\right] . 
\]%
We can know that $\mathbf{t},\mathbf{n},\mathbf{b},\mathbf{e}$ are mutually
orthogonal vector fields satisfying equations

\begin{eqnarray*}
&<&\mathbf{t},\mathbf{t}>_{G}=<\mathbf{n},\mathbf{n}>_{G}=<\mathbf{b},%
\mathbf{b}>_{G}=<\mathbf{e},\mathbf{e}>_{G}=1 \\
&<&\mathbf{t},\mathbf{n}>_{G}=<\mathbf{t},\mathbf{b}>_{G}=<\mathbf{t},%
\mathbf{e}>_{G}=<\mathbf{n},\mathbf{b}>_{G}=<\mathbf{n},\mathbf{e}>_{G}=<%
\mathbf{b},\mathbf{e}>_{G}=0.
\end{eqnarray*}

\bigskip

\newpage 

\textbf{3. MANNHEIM CURVES IN\ GALILEAN SPACE G}$_{4}$\textbf{\bigskip }

In [5], Mannheim curves for Euclidean 4-space are generalized. In this
paper, we have investigated generalization of curves in 4D Galilean space $%
G_{4}$.

\textbf{Definition 3.1. }A special curve $\alpha $ in $G_{4}$ is called a
generalized Mannheim curve if there exists a special Frenet curve $\alpha
^{\ast }$ in $G_{4}$ such that the first normal line at each point of $%
\alpha $ is included in the plane generated by the second and the third
normal line of $\alpha ^{\ast }$ at the corresponding point under $\Psi .$
Here we denote by $\Psi $ a bijection from $\alpha $ to $\alpha ^{\ast }.$
Then the curve $\alpha ^{\ast }$ is called the generalized Mannheim mate
curve of $\alpha .$

A generalized Mannheim mate curve $\alpha ^{\ast }$ is given by the map $%
\alpha ^{\ast }:I^{\ast }\rightarrow G_{4}$ which satisfies the following
equation 
\begin{equation}
\alpha ^{\ast }(t)=\alpha (t)+\gamma (t)\mathbf{n}(t),\quad t\in I. 
\tag{3.1}
\end{equation}%
Here we denote a smooth function on $I$ by $\gamma (t).$ The parameter
should not be an arclength of $\alpha ^{\ast }$. The arclength of $\alpha
^{\ast }$ defined by

\[
t^{\ast }=\int_{0}^{t}\left\Vert \frac{d\alpha ^{\ast }(t)}{dt}\right\Vert
dt 
\]%
Where $t^{\ast }$ is the arclength of $\alpha ^{\ast }.$ For a smooth
function $f:I\rightarrow I^{\ast }$ is given by $f(t)=t^{\ast },$ we have

\[
f^{^{\prime }}(t)=\frac{dt^{\ast }}{dt}=\left\Vert \frac{d\alpha ^{\ast }(t)%
}{dt}\right\Vert =1 
\]%
for $\forall t\in I$. The representation of curve $\alpha ^{\ast }$ with
arclength parameter $t^{\ast }$ is

\[
\alpha ^{\ast }:I^{\ast }\rightarrow G_{4},\quad t^{\ast }\rightarrow \alpha
^{\ast }(t^{\ast }). 
\]%
For the bijection $\Psi :\alpha \rightarrow \alpha ^{\ast }$ defined by $%
\Psi (\alpha (t))=\alpha ^{\ast }(f(t)),$ the reparameterization of $\alpha
^{\ast }$ is given by the following equation

\[
\alpha ^{\ast }(f(t))=\alpha (t)+\gamma (t)\mathbf{n}(t),\quad t\in I 
\]%
where $\gamma (t)$ is a smooth function on $I.$ Then we obtain

\[
\frac{d\alpha ^{\ast }(f(t))}{dt}=\frac{d\alpha ^{\ast }(t^{\ast })}{dt}\mid
_{t^{\ast }=f(t)}=t^{\ast }(f(t)),\quad t\in I. 
\]

\textbf{Theorem 3.1.} If a special Frenet curve $\alpha $ in $G_{4}$ is a
generalized Mannheim curve, the first curvature function $\kappa $ and
second curvature function $\tau $ satisfies the following equation%
\begin{equation}
\kappa (t)=\gamma \tau ^{2}(t)\quad t\in I  \tag{3.2}
\end{equation}%
here we denote a constant number with $\gamma .$

\textbf{Proof.} In the following scheme we show $\alpha $ as a generalized
Mannheim curve and $\alpha ^{\ast }$ as a generalized Mannheim mate curve of 
$\alpha .$%
\[
\begin{array}{ccccc}
&  & \alpha &  & \alpha ^{\ast } \\ 
&  & 
\ddot{}
&  & 
\ddot{}
\\ 
f & : & I & \rightarrow & I^{\ast } \\ 
&  & \downarrow &  & \downarrow \\ 
\Psi & : & G_{4} & \rightarrow & G_{4}%
\end{array}%
. 
\]%
We define a smooth function $f$ by $f(t)=\int \left\Vert \frac{d\alpha
^{\ast }(t)}{dt}\right\Vert dt=t^{\ast }$ is the arclength parameter of $%
\alpha ^{\ast }.$ In addition $\Psi $ is a bijection that is defined by $%
\Psi (\alpha (t))=\alpha ^{\ast }(f(t))$. Then the curve $\alpha ^{\ast }$
is reparameterized as following form%
\begin{equation}
\alpha ^{\ast }(f(t))=\alpha (t)+\gamma (t)\mathbf{n}(t),\quad t\in I 
\tag{3.3}
\end{equation}%
where $\gamma :I\subset \mathbb{R}$ $\rightarrow \mathbb{R}$ is a smooth
function and $\{\mathbf{t},\mathbf{n},\mathbf{b},\mathbf{e}\}$ and $\{%
\mathbf{t}^{\ast },\mathbf{n}^{\ast },\mathbf{b}^{\ast },\mathbf{e}^{\ast
}\} $ are orthogonal vector fields in $G_{4}$ along $\alpha $ and $\alpha
^{\ast },$ respectively.

Differentiating both sides of equation (3.3) with respect to $t,$ we get%
\begin{equation}
\mathbf{t}^{\ast }(f(t))=\mathbf{t}(t)+\gamma ^{^{\prime }}(t)\mathbf{n}%
(t)+\gamma (t)\tau (t)\mathbf{b}(t).  \tag{3.4}
\end{equation}

On the other hand, since the first normal line at the each point of $\alpha $
is lying in the plane generated by the second and the third normal line of $%
\alpha ^{\ast }$ at the corresponding points under bijection $\Psi ,$ the
vector field $\mathbf{n}(t)$ is obtained as follows%
\[
\mathbf{n}(t)=g(t)\mathbf{b}^{\ast }(f(t))+h(t)\mathbf{e}^{\ast }(f(t)) 
\]%
where $g$ and $h$ are some smooth functions on $I\subset \mathbb{R}.$ Taking
into account of the following equation%
\[
<\mathbf{t}^{\ast }(f(t)),g(t)\mathbf{b}^{\ast }(f(t))+h(t)\mathbf{e}^{\ast
}(f(t))>_{G}=0 
\]%
and using (3.4), we have $\gamma ^{^{\prime }}(t)=0.$ Then we decompose the
equation (3.4) as follows%
\begin{equation}
f^{^{\prime }}(t)\mathbf{t}^{\ast }(f(t)=\mathbf{t}(t)+\gamma \tau (t)%
\mathbf{b}(t),  \tag{3.5}
\end{equation}%
that is%
\begin{equation}
\mathbf{t}^{\ast }(f(t)=\mathbf{t}(t)+\gamma \tau (t)\mathbf{b}(t)  \tag{3.6}
\end{equation}%
where $f^{^{\prime }}(t)=1.$

Differentiating both sides of the equation (3.6) with respect to $t\in I$,
we obtain%
\begin{equation}
\kappa ^{\ast }(f(t))\mathbf{n}^{\ast }(f(t))=(\kappa (t)-\gamma \tau
^{2}(t))\mathbf{n}(t)+(\gamma \tau (t))^{^{\prime }}\mathbf{b}(t)+\gamma
\tau (t)\sigma (t)\mathbf{e}(t).  \tag{3.7}
\end{equation}

Using 
\[
<\mathbf{n}^{\ast }(f(t)),g(t)\mathbf{b}^{\ast }(f(t))+h(t)\mathbf{e}^{\ast
}(f(t))>_{G}=0, 
\]%
the coefficient of $\mathbf{n}(t)$ in equation (3.7) vanishes, that is,%
\[
\kappa (t)-\gamma \tau ^{2}(t)=0. 
\]

Then the proof is completed.

\textbf{Theorem 3.2. }Let $\alpha $ be a special Frenet curve such that its
non-constant first and second curvature functions satisfy the following
euation%
\[
\kappa (t)=\gamma \tau ^{2}(t) 
\]%
for all $t\in I\subset \mathbb{R}.$ If the special Frenet curve $\alpha
^{\ast }$ given by the following form%
\[
\alpha ^{\ast }(t)=\alpha (t)+\gamma \mathbf{n}(t) 
\]%
then $\alpha ^{\ast }$ is a generalized Mannheim mate curve of $\alpha .$

\textbf{Proof.} The arclength parameter of $\alpha ^{\ast }$ is given by the
equation

\[
t^{\ast }=\int_{0}^{t}\left\Vert \frac{d\alpha ^{\ast }(t)}{dt}\right\Vert
dt,\quad t\in I. 
\]%
Let us assume that%
\[
\kappa (t)=\gamma \tau ^{2}(t), 
\]%
then we obtain $f^{^{\prime }}(t)=1$, $t\in I.$

Differentiating the equation $\alpha ^{\ast }(f(t))=\alpha (t)+\gamma 
\mathbf{n}(t)$ with respect to $t$ the we get%
\[
f^{^{\prime }}(t)\mathbf{t}^{\ast }(f(t)=\mathbf{t}(t)+\gamma \tau (t)%
\mathbf{b}(t). 
\]%
Then we can see%
\begin{equation}
\mathbf{t}^{\ast }(f(t)=\mathbf{t}(t)+\gamma \tau (t)\mathbf{b}(t),\quad
t\in I.  \tag{3.8}
\end{equation}

\bigskip

Differentiating the last equation with respect to $t$ is%
\begin{equation}
\kappa ^{\ast }(f(t))\mathbf{n}^{\ast }(f(t))=(\kappa (t)-\gamma \tau
^{2}(t))\mathbf{n}(t)+(\gamma \tau (t))^{^{\prime }}\mathbf{b}(t)+\gamma
\tau (t)\sigma (t)\mathbf{e}(t).  \tag{3.9}
\end{equation}

From the assumption, we obtain%
\[
\kappa (t)-\gamma \tau ^{2}(t)=0. 
\]

Then, the coefficient of $\mathbf{n}(t)$ in the equation (3.9) is zero. One
can see from the equation (3.8) that $\mathbf{t}^{\ast }(f(t)$ is a linear
combination of $\mathbf{t}(t)$ and $\mathbf{b}(t).$ In addition, from
equation (3.9), $\mathbf{n}^{\ast }(f(t))$ is given by linear combination of 
$\mathbf{b}(t)$ and $\mathbf{e}(t).$ On the otherhand, $\alpha ^{\ast }$ is
a special Frenet curve that the vector $\mathbf{n}(t)$ which satisfies the
following linear combination of $\mathbf{t}^{\ast }(f(t)$ and $\mathbf{n}%
^{\ast }(f(t)).$

Therefore, the first normal line $\alpha $ lies in the plane generated by
the second and third normal line of $\alpha ^{\ast }$ at the corresponding
points under the $\Psi $ bijection which is defined by \newline
$\Psi (f(t))=\alpha ^{\ast }(f(t)).$ The proof is completed.

\textbf{Remark 3.1.} In 4D Galilean space $G_{4},$ a special Frenet curve $%
\alpha $ with curvature functions $\kappa $ and $\tau $ satisfying $\kappa
(t)=\gamma \tau ^{2}(t),$ it is not clear that a smooth curve $\alpha ^{\ast
}$ given by (3.1) is a special Frenet curve. The reverse of Theorem 3.1 is
still a great puzzle for the authors.

\bigskip

\newpage 

\textbf{References}

[1] H. Balgetir, M. Bekta\c{s} and J. Inoguchi, \emph{Null Bertrand curves
in Minkowski 3-space and their characterizations,}, Note Math., 23 no. 1,
7-13, 2004.

[2] FN. Ekmekci and K. \.{I}larslan,, \emph{\ On Bertrand curves and their
characterization}, Differ. Geom. Dyn. Syst.(electronic), vol.3, no.2, 2001.

[3] R. Blum \emph{\ A remarkable class of Mannheim curves,} Canad. Math.
Bull., 9, 223-228, 1966.

[4] H. Liu and F. Wang, \emph{Mannheim Partner curves in 3-space}, Journal
of Geometry, 88, 120-126, 2008.

[5] H. Matsuda and S. Yorozu \emph{On generalized Mannheim curves in
Euclidean 4-space}, (English) Nihonkai Math. J., 20, no. 1, 33-56, 2009.

[6] S.Yilmaz \emph{Construction of the Frenet-Serret frame of a curve in 4D
Galilean space and some applications}, Int. Journal of the Physical Sciences
Vol. 5(8), pp. 1284-1289, 2010.

[7] O. Roschel,\emph{Die geometrie Des Galileischen Raumes} Berichte der
Math.-Stat. Sektion im Forschungszentrum Graz Ber., 256: 1-20, 1986.

[8] IM. Yaglom, \emph{A Simple Non-Euclidean Geometry and Its Physical Basis}%
, Springer-Verlag, New York, 1979.

[9] A.O. \"{O}\u{g}renmis, H.\"{O}ztekin, M. Erg\"{u}t, \emph{Bertrand
curves in Galilean space and their characterizations}, Kragujevac J. Math.,
32: 139-147, 2009.

[10] A.O. \"{O}\u{g}renmis, M. Erg\"{u}t, M. Bekta\c{s}, \emph{On the
Helices in the Galilean Space G3}, Iran. J. Sci. Tech. Trans. A Sci., 31(2):
177-181, 2007.

[11] S.Ersoy, M.Akyi\u{g}it and M.Tosun, \emph{\ A Note on Admissible
Mannheim Curves in Galilean Space G3}, International J.Math. Combin. Vol.1,
88-93, 2011.

\end{document}